\documentclass[twoside,leqno]{article}

\usepackage[letterpaper]{geometry}

\usepackage{siamproceedings}

\usepackage[T1]{fontenc}
\usepackage{amsfonts}
\usepackage{graphicx}
\usepackage{epstopdf}
\usepackage{enumitem}
\usepackage{algorithmic}

\ifpdf
  \DeclareGraphicsExtensions{.eps,.pdf,.png,.jpg}
\else
  \DeclareGraphicsExtensions{.eps}
\fi

\newsiamremark{remark}{Remark}
\newsiamremark{hypothesis}{Hypothesis}
\crefname{hypothesis}{Hypothesis}{Hypotheses}
\newsiamthm{claim}{Claim}

\usepackage{amsopn}

\usepackage{bbm}
\usepackage{xcolor}

\newcommand\barc{A}

\newcommand\dd{\displaystyle}
\newcommand\dv{\mathop{\rm div}}

\newcommand{\m}[1]{\mathbbm{#1}}

\newcommand{\q}[1]{\mathcal{#1}}
\newcommand{\ps}{{\partial_s}}

\newcommand\RR{{\cal R}}
\renewcommand\SS{{\cal S}}

\newcommand{\xx}{{a}}

\begin{document}



%
\newcommand\relatedversion{}
\renewcommand\relatedversion{\thanks{The full version of the paper can be accessed at \protect\url{https://arxiv.org/abs/0000.00000}}} 

\title{\Large \textbf{Soliton resolution at blow-up for the subconformal nonlinear
    wave equation}}
\author{Hatem Zaag \thanks{Universit\'e Sorbonne Paris Nord,  LAGA,
    CNRS (UMR 7539), F-93430, Villetaneuse, France (\email{hatem.zaag@univ-paris13.fr})}} 

\date{}

\maketitle

\fancyfoot[R]{\scriptsize{Copyright \textcopyright\ 20XX by SIAM\\
Unauthorized reproduction of this article is prohibited}}

\begin{center}
\textit{Dedicated to the memory of the late Professor Robert V. Kohn.}
\end{center}

\begin{abstract} 
  We consider the nonlinear wave equation (NLW) with a superlinear pure
  power nonlinearity in the subconformal and conformal cases.
  Under some conditions
  on initial data, this equation is known to have solutions which blow up in
  finite time. Two questions are then relevant:
  (i) the classification of all possible blow-up
  behaviors; (ii) the construction of examples of blow-up solutions.
  As we will show, the situation is entirely settled in the
  one-dimensional case, where we fully solve the famous \textit{``Soliton
  Resolution Conjecture''}. Various extensions to higher dimensions and
  to perturbed versions of (NLW) are given, including the
  construction of a solution in 2-d, with a nearly pyramidal blow-up
  graph. 
  Carrying out this program was made possible thanks to a synergy of
  techniques from the PDE theory, Mathematical Physics and Analysis, including
  ODE techniques, spectral theory and energy methods.
  All over this presentation, we will insist on connections with the
  study of other types of PDEs, in particular in the parabolic
  case. Surprisingly enough, in spite of the difference between
  parabolic and hyperbolic equations at the linear level, the
  nonlinear nature
  brings in a strong unity - both in the results and in the
  methods - between these two important classes of PDEs.
%
%
\end{abstract}


\medskip



\section{Introduction.}

\subsection{Blow-up in nature.}

Differential Equations are commonly used to model various real-world
phenomena. Under some circumstances, a physical quantity may
undergo a dramatic change of magnitude, in some finite time. This is
for instance the case with chemotaxis, a situation where the amoebas
of the species Dictyostellium Discoideum 
attract each other by means of some chemical they secrete,
when ressources become rare (see Horstmann
\cite{Hjdmv03}). In that case, their concentration becomes very high
in some localized region, in a relatively short time. Mathematically,
we say that the concentration ``\textit{blows up in finite
  time}''.
Finite time blow-up
occurs in a large variety of situations,
including micro-electronics, with the ''touch-down'' phenomenon in
Micro-Electro-Mechanical Systems or MEMS (see Kavallaris and Suzuki
\cite{KSspringer18}). Many situations from nonlinear optics to phase
transitions
modeled by the Complex Ginzburg-Landau equation (CGL) may exhibit
blow-up as well (see Popp \textit{et al} \cite{PSKKphysd98}). We
also encounter blow-up in general relativity (see Donninger, Schlag and
Soffer \cite{DSScmp12}), as well as in self-focusing waves in nonlinear
optics (see Bizo\'n, Chmaj and Szpak \cite{BCSjmp11}).

\medskip

Many of the
above-mentioned problems are difficult to handle mathematically,
because they lack ``structure'' (conserved or dissipated energy
functionals, comparison principles, invariance under some
transformations, nice spectral properties). For
that reason, a large focus has been first made on blow-up for model
problems, which capture essential features of physically relevant
equations, while being simple enough to be accessible to mathematical
analysis. Among them, we cite the Nonlinear Wave equation (NLW) which is
the main object of our paper. We will also mention the blow-up
analysis for the semilinear heat equation, which was quite inspiring
for the study of NLW. Addressing these equations is already a
difficult task, as one can see in this paper.

\medskip

Some of our techniques and results were later successfully adapted to physically
relevant problems (see for example our contributions on CGL in
\cite{DNZarma24}, the MEMS model in \cite{DZm3as19} and also the
Keller-Segel model for chemotaxis in \cite{NNZapde25}).

\subsection{Blow-up in Differential Equations.}

Under some conditions, nonlinear Differential Equations
may exhibit blow-up in finite time. For example, if one considers the
following Ordinary Differential Equation (ODE)
 \begin{equation}\label{ode}
u'= f(u),
 \end{equation}
 where $u= u(t) \in \m R$, $t\ge 0$
 and $f:\m R \to \m R$ is positive and continuous, it is easy to see that whenever
 \[
\int_1^\infty \frac 1 f <+\infty,
 \]
 then, equation \eqref{ode} has a blow-up solution $u(t)$, defined for all
 $t\in [0,T)$ for some finite $T>0$, such that
 \[
u(t) \to \infty\mbox{ as } t\to T.
 \]
 Partial Differential Equations (PDE) may have blow-up solutions as
 well. This is the case for the following Nonlinear Wave (NLW) equation  
\begin{equation}\label{eqnlw}
\partial_t^2 u = \Delta u +|u|^{p-1}u,
\end{equation}
where $u=u(x,t)$, $x\in \m R^N$ and $t\ge 0$. The exponent $p$ will
be considered in the subconformal range: 
\begin{equation}\label{condpwave}
  1<p<p_c
\end{equation}
and in some statements in the conformal case too:
\begin{equation}\label{condpwave'}
  N\ge 2 \mbox{ and }
  p =p_c,
  \end{equation}
  where
  \[
p_c = \frac{N+3}{N-1} \mbox{ if } N\ge 2\mbox{ and } p_c = +\infty
\mbox{ if } N=1.
  \]
  The exponent $p_c$ is referred to as the ``conformal'' exponent,
because (NLW) is invariant under the so-called ``conformal''
transformation precisely when $p=p_c$ and only for that exponent
(see Merle and Zaag \cite{MZma05}). Note
that whenever $N\ge 3$, it holds that
$p_c <p_s$,
the Sobolev critical exponent defined by
\[
p_S = \frac{N+2}{N-2} \mbox{ if } N\ge 3\mbox{ and } p_S = +\infty
\mbox{ if } N\le 2.
  \]
From Levine \cite{Ltams74}, we know that the solution of \eqref{eqnlw}
blows up (in a sense that will be given in Section \ref{secwave}), if the following energy
functional
\[
  \dd\int_{{\m R}^N}\left(\frac 12 |\partial_t u(x,t)|^2+\frac 12|\nabla u(x,t)|^2
    -\frac 1{p+1}|u(x,t)|^{p+1}\right)dx
\]
(which is conserved) is negative. Equipping equation \eqref{eqnlw}
with initial data $(u, \partial_tu)(x,0) = (\lambda \varphi, \psi)(x)$
for some regular compactly supported $\varphi\not\equiv 0$ and $\psi$ and large
enough $\lambda>0$, one can produce a solution with negative energy,
which blows up.

\medskip

Equation \eqref{eqnlw} enjoys a nice physical property: the finite
speed of propagation, in the sense that the value of $u$ at some
$(x_0,t_0)$ depends only on the values of $u$ in the backward light
cone $\q C_{x_0,t_0,1}$, where for any $x\in \m R^N$, $t\ge 0$ and
$\delta>0$, we introduce the cone
\begin{equation}\label{defcone}
{\cal C}_{x,t,\delta}=\{(\xi,\tau)\;|\; 0\le \tau< t- \delta|\xi-x|\}.
\end{equation}
In particular, any change in initial data outside the section of $\q
C_{x_0,t_0,1}$ at $\tau=0$ makes no change in the value of $u$ at its
vertex $(x_0,t_0)$. This property will prove to be essential in our
blow-up analysis.
In particular, it enables us, through a simple cut-off technique, to
derive space dependent blow-up solutions for \eqref{eqnlw} from
blow-up solutions of the associated ODE $u''=u^p$. 

 \medskip

 In the following, we aim at addressing the blow-up question for \eqref{eqnlw} in the
 range \eqref{condpwave}-\eqref{condpwave'}. This will lead us to the so called
 \textit{``Soliton Resolution Conjecture''} (or SR conjecture for short), which we completely solve
 for $N=1$ (see Subsection \ref{secasymp}). We will also
provide a
 blow-up solution in  higher dimensions, which obeys that
 conjecture, showing a nearly pyramidal shape for the boundary of the domain
 of definition (see Subsection \ref{secpyr}).

 \medskip

Interestingly, our strategy is largely inspired by the parabolic
 case, namely the case of the following semilinear heat equation:
 \begin{equation}\label{eqchaleur}
\partial_t u = \Delta u +|u|^{p-1}u,
\end{equation}
where $u=u(x,t)$, $x\in \m R^N$, $t\ge 0$, in the Sobolev subcritical
range 
\begin{equation}\label{condpchaleur}
  1<p<p_S.
\end{equation}
This inspiration may seem surprising, given the large differences in
the behavior between the linear versions of both equations.
It is as if the nonlinear nature brings in a
strong unity in the results and in the methods, between
these two equations.
For that reason, we will briefly review the main blow-up results
for equation \eqref{eqchaleur}, before addressing the case of equation
\eqref{eqnlw}.

\medskip

Before all that, we will first review the relevant questions about
blow-up, then introduce the SR conjecture.

\subsection{Relevant questions about blow-up.} \label{secquest}
Two general questions arise when dealing with the behavior of PDEs, in
particular at blow-up:

\medskip

- {\bf Classification}:
Given a PDE, is it possible to provide a classification of all the
possible behaviors it may exhibit? 
Often, due to the difficulty of such a general classification, one may
instead focus on subclasses of solutions (for instance, global solutions,
or solutions that blow up in finite time). One may also further
restrict the study, for example to radial solutions, or by seeking
such a classification for particular initial data, notably in the
neighborhood of a universal object, such as a “soliton” (see
Subsection \ref{secsrc} for a definition).

\medskip

- {\bf Construction of examples}: Given a PDE, can one
construct examples of solutions exhibiting a particular behavior, or, in the usual terminology, ``{\it with a prescribed behavior}''?\\
If the answer to such a question serves to confirm, through examples,
a general classification (if it exists), it is precisely in the
absence of such a classification that the construction proves to be
most crucial. Indeed, in such a case, the constructed examples often
constitute the only known examples of solutions.

\subsection{The Soliton Resolution Conjecture.} \label{secsrc}
Whether adopting the classifying or the constructive point of
view, asymptotic behavior of solutions to PDEs
is an essential feature.
Even at blow-up, thanks to some transformations (see
below \eqref{defwondes} for the case of (NLW)), we may reduce the
question of blow-up asymptotics to long-time asymptotics for some
related PDE (as in \eqref{eqw-ondes0}).

\medskip

As far as the asymptotic behavior is concerned for dispersive equations
(NLW, NLS, KdV\footnote{NLS stands for the Nonlinear Schr\"odinger
  equation, and KdV for the Korteweg - de Vries equation.},...), 
an important
conjecture has been circulating in the community: the
\textit{``Soliton Resolution Conjecture''}. According to Terence Tao
on his blog \cite{Tblog08},
\textit{``this conjecture (...) asserts, roughly speaking, that any
  reasonable (e.g. bounded energy) solution to such equations
  eventually resolves into a superposition of a radiation component
  (which behaves like a solution to the linear (...) equation) plus a
  finite number of “nonlinear bound states” or “solitons”''.}

\medskip

Solitons are indeed encountered in various dispersive equations. We
may define a soliton as a stable solitary wave (i.e. a function of the
form $Q(x-ct)$), which preserves its
shape (up to small perturbations) after a collision with another
solitary wave. Historically, solitary waves were first discovered
by Korteweg and de Vries \cite{KdVpm95} in 1895 in the context of
fluid dynamics.

\medskip

The SR conjecture has attracted a lot of attention in the last twenty
years, and no list of authors can be exhaustive.
The first results
were obtained in
perturbative cases, i.e. for initial data close to solitons. More
recently, further results have been established for genuinely
non-perturbative cases (see Duyckaerts, Kenig and Merle
\cite{DKMjems12} and \cite{DKMcjm13} for the energy-critical Nonlinear
Wave equation \eqref{eqnlw}, and also their work together with Jia 
\cite{DJKMimrn18} for the energy-critical wave maps equation).

\medskip

For parabolic equations, it is possible to formulate an analogue of the
SR conjecture, even though the notion of  ``soliton'' may appear somewhat
inappropriate, outside the field of dispersive
equations. Nevertheless, universal objects, namely solutions of
stationary equations associated with the PDE, often arise as natural
candidates for the  ``profile'' in parabolic equations, in particular
near blow-up (see Bressan \cite{Biumj90} and \cite{Bjde92},
Bricmont-Kupiainen \cite{BKnonl94}, Herrero-Velázquez \cite{HVihp93}
for the semilinear heat equation).  

\section{A brief sketch of blow-up for the semilinear heat
  equation.}\label{secchaleur}

Here, we consider the semilinear heat equation \eqref{eqchaleur} under
the condition \eqref{condpchaleur}. We consider initial data in
$L^\infty(\m R^N)$, where the local-in-time Cauchy problem follows
from a simple fixed point argument. Accordingly, we say that a
solution blows up when its $L^\infty$ norm goes to infinity. The first
sufficient conditions 
for blow-up go back to Kaplan \cite{Kcpam63}, Fujita \cite{Fjfsut66}
and Levine \cite{Larma73}. The literature about equation
\eqref{eqchaleur} is huge, and no list can be exhaustive. We refer the
interested reader to the book by Quittner and Souplet \cite{QSbirk19}
and the references therein. In the following, we will review some
blow-up results from the eighties and the nineties, which were
inspiring for us, while addressing blow-up for the Nonlinear Wave
equation \eqref{eqnlw}.

\subsection{Classification results.}
The first classification results for \textit{general} blow-up solutions
to equation \eqref{eqchaleur} were achieved by Giga and Kohn in the
eighties in \cite{GKcpam85}, \cite{GKiumj87} and \cite{GKcpam89},
under a positivity condition or some further constraint on $p$ inside
the Sobolev subcritical range \eqref{condpchaleur}. Later in 2004,
Giga, Matsui and Sasayama \cite{GMSiumj04} were able to derive the blow-up rate (see
\eqref{boundw} below) with no restriction. As a consequence, all the
program carried out by Giga and Kohn became immediately valid in the
whole range shown in  \eqref{condpchaleur}. 

Given a solution $u(x,t)$ to \eqref{eqchaleur} which blows up at time
$T>0$, Giga and Kohn use the so-called \textit{similarity variables}
transformation, which was first introduced by Hocking and Stewartson \cite{HSprsl72},
in the more general context of the Complex Ginzburg-Landau (CGL) equation:
\begin{equation}\label{defwchaleur} 
w_{x_0}(y,s) = (T-t)^{\frac 1{p-1}}u(x,t),\;\;y=\frac {x-x_0}{\sqrt{T-t}},\;\;s=-\log(T-t).
\end{equation}
In this setting, the function $w_{x_0}$ (or $w$ for short) satisfies
the following equation, for all $s\ge s_0\equiv- \log T$ and $y\in \m R^N$:
\begin{equation}\label{eqwchaleur}
\partial_s w = \frac 1 \rho \dv(\rho \nabla w) - \frac w{p-1}+|w|^{p-1}w,
\end{equation}
where
\[
\rho(y) = e^{-\frac{|y|^2}4}.
\]
The first major finding of Giga and Kohn was to notice the existence
of a Lyapunov functional for $w$, defined by
\begin{equation}\label{defE}
E_h(w) = \dd \int_{\m R^N} \left( \frac 12 |\nabla w|^2 +
  \frac{w^2}{2(p-1)} - \frac{|w|^{p+1}}{p+1}\right)\rho dy,
\end{equation}
in the sense that
\begin{equation}\label{dissipE}
\frac d{ds} E_h(w(s)) = - \int_{\m R^N}(\partial_s w)^2 \rho dy.
\end{equation}
Using some blow-up criterion, they showed that
\begin{equation}\label{boundE0}
  \forall s\ge s_0,\;\; 0\le E_h(w(s)) \le E_h(w(s_0)).
\end{equation}
Thanks to energy techniques and interpolation, Giga and Kohn together
with Giga, Matsui and Sasayama \cite{GMSiumj04}
showed that
\begin{equation}\label{boundw}
\forall s\ge s_0,\;\;\|w(s)\|_{L^\infty} \le C.
\end{equation}
From \eqref{dissipE} and \eqref{boundE0}, 
the Giga-Kohn program continues with
the convergence of $E_h(w(s))$ as $s\to \infty$,
which reasonably indicates that $\partial_s w\to 0$ as $s\to
\infty$, thanks to \eqref{dissipE}. This implies the convergence of
$w(y,s)$ to some stationary solution of \eqref{eqwchaleur}. Then,
thanks to some energy and Pohozaev type identities, they showed
that in the Sobolev subcritical range \eqref{condpchaleur}, the only
stationary solutions are given by $0$, $\kappa$ and $-\kappa$, where
\[
\kappa=(p-1)^{-\frac 1{p-1}}.
\]
Later, Herrero and Vel\`azquez \cite{HVcpde92} together with Filippas
and Kohn \cite{FKcpam92} derived the \textit{``blow-up profile''} for
$N=1$, in the sense that one of the following scenarios happens (up to
replacing $u$ by $-u$):  
\begin{itemize}
\item[-] either  $u(x,t) \equiv \kappa (T-t)^{-\frac 1{p-1}}$,
\item [-] or $u(x,t) \sim (T-t)^{-\frac 1{p-1}} \left(p-1 +
    \frac{(p-1)^2}{4p}\frac{|x|^2}{(T-t)|\log(T-t)|}\right)^{-\frac
    1{p-1}}$ as $t\to T$,
  \item [-] or $u(x,t) \sim (T-t)^{-\frac 1{p-1}}
  \left(p-1 + b\frac{|x|^{m}}{(T-t)}\right)^{-\frac
    1{p-1}}$ as $t\to T$,
\end{itemize}
for some $b>0$ and some even integer $m\ge 4$. The functions shown on
the right-hand sides of the two latter statements are referred to as
\textit{``blow-up profiles''}. This statement was later generalized in
higher dimensions by Vel\'azquez in \cite{Vcpde92} (see also Filippas
and Liu \cite{FLihp93}). 

\subsection{Construction results.}

Following the classification we have just given, a natural question
arises: do we have examples of initial data leading to those
scenarios?

\medskip

When $N=1$, the answer came from Bricmont and Kupiainen
\cite{BKnonl94} (see also Herrero and Vel\'azquez \cite{HVihp93} as
well as Merle and Zaag \cite{MZdmj97}).
Roughly speaking, the
method in \cite{BKnonl94} and \cite{MZdmj97} follows a two-stage
strategy (see Tayachi and Zaag \cite{TZtams19}
where the strategy is
implemented
for a double source nonlinear heat
equation): 
\begin{itemize}
\item[-] A formal approach, relying on approximations of the PDE in various
regions of space, together with a matching asymptotics step. Its
outcome is the determination of the candidate for the blow-up profile.
\item[-] A rigorous approach, where one linearizes the equation around the
  candidate of the blow-up profile. Using the spectral properties of
  the linearized operator, one sees that the projection on the negative
  part of the spectrum can be easily controlled, from the spectral
  properties. This way, only the 
  projection on the nonnegative part of the spectrum remains to be
  controlled. That part is in fact finite dimensional. Thus, one only
  needs to deal with a finite dimensional problem, whose solution is
  easily derived from degree theory.
\end{itemize}
The reduction to finite dimensions yields the so called
\textit{``co-dimensional stability''} of the considered profile. Note
that the rigorous method is in fact linked to the \textit{center
  manifold theory}, which (in general) cannot be directly applied,
mainly because one cannot find a suitable functional space, where the
linear operator is easily understood, while the nonlinear term is
continuous in that space (see Filippas and Kohn \cite{FKcpam92} for a
more detailed discussion). However, in other contexts, in particular
with supercritical problems, some authors were able to
introduce some spaces where the linear operator was well understood
while the nonlinear part is continuous. For examples in the parabolic
case, see Biernat, Donninger and Sch\"orkhuber \cite{BDScvpde17},
Golgi\'c and Sch\"orkhuber \cite{GScpde20}, \cite{GSarma24}, Golgi\'c,
Kistner and Sch\"orkhuber \cite{GKScvpde24}, Ghoul, Ibrahim and Nguyen
\cite{GINapde19}. 

\bigskip

In higher dimensions $N\ge 2$, a notion of degeneracy appears in the
predicted blow-up profiles for equation \eqref{eqchaleur}. When the
predicted profile is not degenerate, we 
know from Amadori \cite{Adie95} examples of initial data leading to
that behavior. In the degenerate case, the question was recently
solved in a joint work with Merle in \cite{MZjems24}. Handling degenerate
cases required a very fine technique of integration in time for the PDE,
component by component. It leads in particular to the construction of
a solution $u(x,t)$ to equation \eqref{eqchaleur} in 2 space
dimensions which blows up at some time $T>0$ only at the origin, with a
\textit{cross-shaped} blow-up profile, in the sense that
\[
u(x,T) \sim\left[\frac{(p-1)^2}\kappa (x_1^2x_2^2+\delta(x_1^6+x_2^6))\right]^{-\frac 1{p-1}}\mbox{ as }x\to 0,
\]
where for all $x\neq 0$, $u(x,T) = \dd\lim_{t\to T}u(x,t)$.

\section{The Nonlinear Wave Equation.}  \label{secwave}

We consider equation \eqref{eqnlw} in the subconformal and conformal
ranges \eqref{condpwave}-\eqref{condpwave'}, with initial data  $(u, \partial_t u)(x,0)$
in the space $H^1\times L^2(\m R^N)$.
Some generalizations of this equation will be considered as well.

\medskip

We first review some nice properties for equation \eqref{eqnlw},
which will impact our blow-up analysis.

\medskip

First, we recall the finite speed of propagation, already mentioned in
the introduction, which makes the analysis in backward light cones
relevant (take $\delta=1$ in \eqref{defcone}). Then, we cite the
invariance by translation, both in space and in time. The latter is
responsible of a positive eigenvalue ($\lambda=1$) in the spectrum
of linearized operators. Finally, we introduce the following Lorentz
transform, relevant in one space dimension, which keeps the equation
invariant as well:

\medskip

For $N=1$ and for any $d\in(-1,1)$, the function $U$ defined by 
\[
U(x',t')=u(x,t) \mbox{ where }x'= \frac{x+dt}{\sqrt{1-d^2}}\mbox{ and }t'=\frac{t+dx}{\sqrt{1-d^2}}
\]
is also a solution of \eqref{eqnlw}. This invariance is responsible
for a null eigenvalue in the spectrum, difficult to handle. Some
modulation technique will be needed to overcome that difficulty.

\medskip

The Cauchy problem for equation \eqref{eqnlw} can be solved in the
space ${\rm H}^1\times {\rm L}^2$ (see for instance Ginibre, Soffer
and Velo \cite{GSVjfa92} and the book \cite{STnyu98} by Shatah and
Struwe).  If the solution is not global in time, then we call it a
blow-up solution. As already mentioned in the introduction, the
existence of blow-up solutions follows from Levine's criterion in
\cite{Ltams74} or from localizing ODE solutions.
More blow-up results can be found in
Caffarelli and Friedman \cite{CFtams86}, \cite{CFarma85},
Alinhac \cite{Apndeta95}
and Kichenassamy and Littman \cite{KL1cpde93}, \cite{KL2cpde93}.

\medskip

As we have already mentioned in the introduction, equation \eqref{eqnlw}
can be considered as a lab model for blow-up in 
hyperbolic equations, because it captures features common to a whole
range of blow-up problems arising in various nonlinear physical
models, in particular in general relativity (see Donninger, Schlag and
Soffer \cite{DSScmp12}), and also for self-focusing waves in nonlinear
optics (see Bizo\'n, Chmaj and Szpak \cite{BCSjmp11}).

\medskip

Consider $u(x,t)$ a blow-up solution to equation \eqref{eqnlw}.
From Alinhac \cite{Apndeta95},
one can define a 1-Lipschitz continuous \textit{blow-up graph} $\Gamma=\{(x,T(x))\}$
such that the domain of definition of $u$ (also called \textit{the maximal influence
domain} of $u$) is given by 
\begin{equation}\label{defdu}
D_u=\{(x,t)\;|\; 0\le t< T(x)\}.
\end{equation}
This definition shows that the solution is regular up to $\bar
T=\dd\inf_{x\in {\m R^N}}T(x)$, which can be referred to as the
``first'' blow-up time.
\begin{figure}[h]\label{figDu}
\begin{center}
\includegraphics[width=0.45\textwidth]{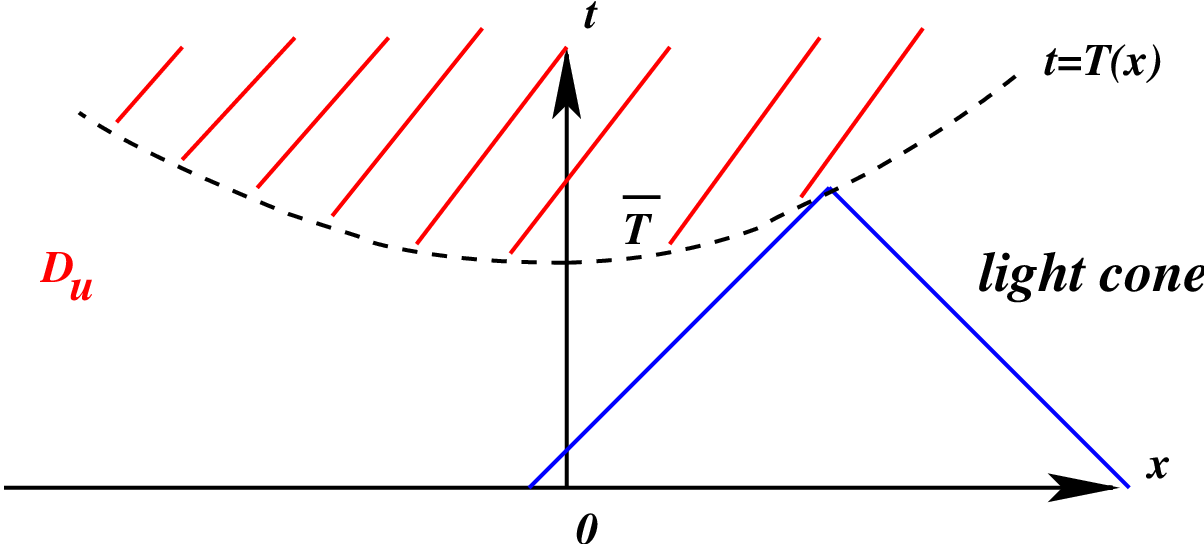}
\caption{Domain of definition
 of the Nonlinear Wave equation}
\end{center}
\end{figure}
Whenever the set
$\bar E \equiv \{x \in \m R^N\;|\;\bar T(x) =\bar T\}$ is not empty,
for example for compactly supported initial data, its 
elements are the first points where the solution blows up. Unlike the
situation of the semilinear heat equation \eqref{eqchaleur}, where we
have only one blow-up time, and where the solution cannot be extended
beyond that time (see Baras and Cohen \cite{BCjfa87}), the solution of
the Nonlinear Wave equation \eqref{eqnlw} may continue to exist beyond $\bar
T$,
in regions of space away from $\bar
E$. Indeed, for points $x\not\in\bar E$, the solution will blow up at
some later time, namely $T(x)$, which may be referred to as a
``local'' blow-up time. In other words, the solution blows up
everywhere in space, however, at different blow-up times (see Figure
\ref{figDu}). 
 
\medskip

The $1$-Lipschitz character of $\Gamma$ is a consequence of the finite speed of
propagation. Indeed, by construction, $D_u$ can be seen as the union
of all backward light cones ${\cal C}_{a,T(a),1}$ \eqref{defcone} for
any $a\in {\m R}^N$, which implies that $x\mapsto T(x)$ is Lipschitz
with a constant bounded by the common slope of the cones, namely $1$.  

\medskip

As $1$ appears to be the maximal possible value for the local
Lipschitz constant of $x\mapsto T(x)$, one may expect some kind of
tangency between the blow-up graph and some backward light cone, which
may affect the blow-up analysis. For that reason, we distinguais two
categories of points in $\m R^N$, according to the value of the local
Lipschitz constant.

\medskip

\begin{figure}[htbp]\label{Fignc}
  \centering
  \includegraphics[width=0.5\columnwidth]{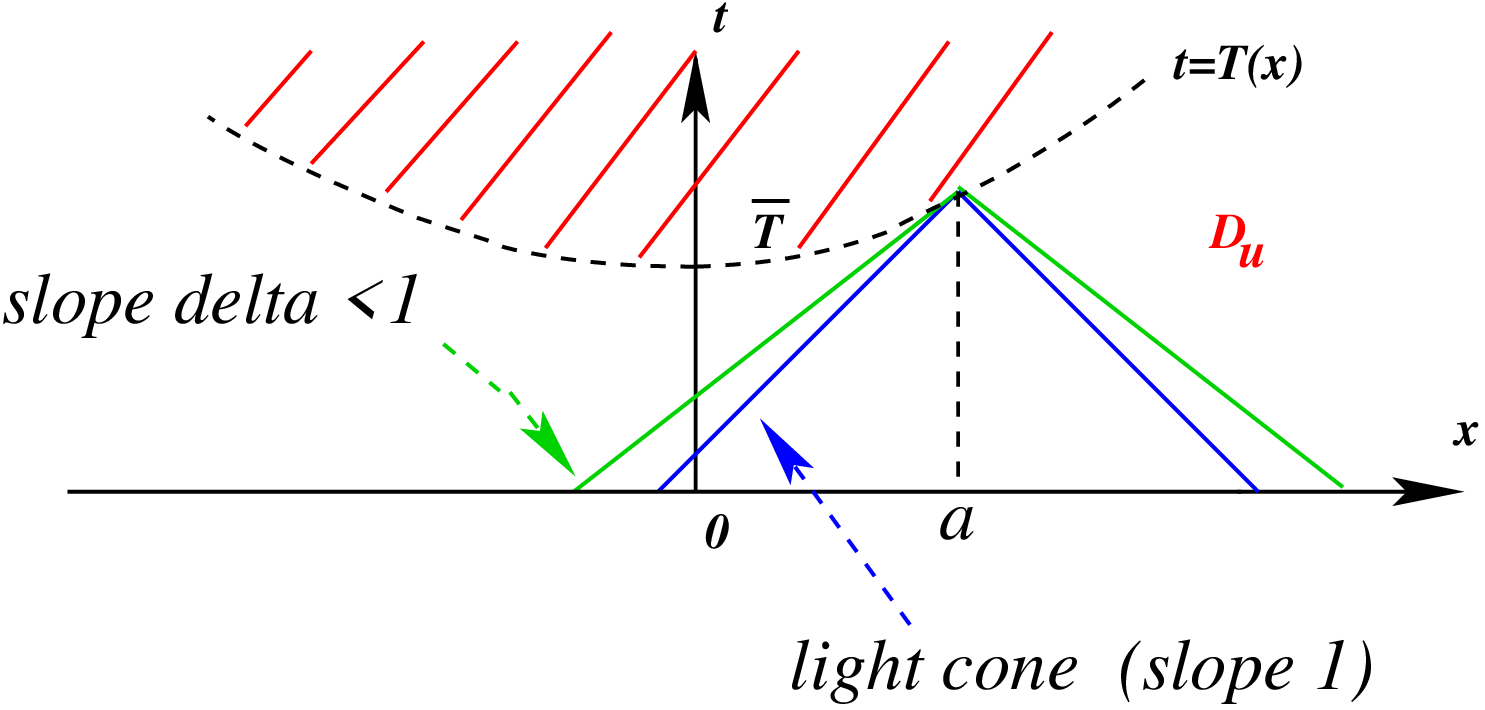}
  \caption{The point $a$ is non-characteristic.}
\end{figure}

More precisely, a point $a\in \m R^N$ is called non-characteristic (or {\it regular})
if there exists $\delta\in [0,1)$
such that $u$ is defined on ${\cal C}_{a, T(a),\delta}$ (see Figure \ref{Fignc}).
If not, then we call $a$ a characteristic point (or a {\it singular}
point).
We denote by $\RR$ (resp. $\SS$) the set of non-characteristic (resp. characteristic) points. Note then that 
\[
\RR \cup \SS = {\m R}^N.
\]
We will show subsequently that this distinction between $\RR$ and
$\SS$ is relevant,
since the asymptotic behavior is different
according to whether we are near one set or the other (see Subsection
\ref{secasymp} below).

\medskip

We remark that $\RR \not =\emptyset$, for any blow-up
solution. Indeed, this is clear for compactly supported initial data, since one
easily sees that $x\mapsto T(x)$ has a minimum, which happens to be a
non-characteristic point (with $\delta=0$). For more general data, see
Proposition 3.5 in \cite{MZcmp08}.

\medskip

Accordingly, four fundamental questions arise as far as blow-up is
concerned for equation \eqref{eqnlw}:

\medskip

\noindent
- \textbf{Question 1}: \textit{Can we estimate the ``size'' (or
  blow-up rate) of the solution at blow-up?}\\
- \textbf{Question 2}: \textit{Can we derive the asymptotic behavior of solutions at blow-up?}\\
- \textbf{Question 3}: \textit{Can we derive more regularity for the blow-up graph?}\\
- \textbf{Question 4}: \textit{Are there blow-up solutions
  with $\SS\not =\emptyset$?}

\medskip

In a series of papers with Merle (\cite{MZajm03}, \cite{MZma05}, \cite{MZimrn05},
\cite{MZjfa07}, \cite{MZcmp08}, \cite{MZajm12}, \cite{MZdmj12}) and
with C\^ote \cite{CZcpam13}, we
made several contributions concerning these issues, in particular 
the description of the blow-up graph and the blow-up behavior in
similarity variables, bringing the knowledge of  the blow-up problem
for equation \eqref{eqnlw} to the level reached for the semilinear
heat equation by Giga and Kohn \cite{GKcpam89}, Chen and Matano
\cite{CMjde89}, Herrero and Vel\'azquez \cite{HVihp93}. In particular,
we completely solve the \textit{``Soliton Resolution Conjecture''} for equation
\eqref{eqnlw} when $N=1$.

\medskip

The following subsections are devoted to the above-mentioned
questions. The last one
is devoted to an example of a solution in 2
space dimensions, with a nearly pyramidal shape for the blow-up
graph.

\subsection{A Lyapunov functional and blow-up rate.}

A fundamental feature of equation \eqref{eqnlw} lays in the
existence of a Lyapunov functional in the so-called \textit{similarity
  variables} setting defined for each $x_0\in \m R^N$ as follows:
\begin{equation}\label{defwondes}
w_{x_0}(y,s) = (T(x_0)-t)^{\frac 2{p-1}}u(x,t),\;\;y=\frac {x-x_0}{T(x_0)-t},\;\;s=-\log(T(x_0)-t).
\end{equation}
This is in fact the hyperbolic version of the transformation
\eqref{defwchaleur} first introduced in the parabolic setting.
After the transformation \eqref{defwondes}, we recover the following
equation satisfied by $w_{x_0}$ (we write $w$ for simplicity):
\begin{equation}\label{eqw-ondes0}
\partial^2_{s}w= {\q L} w-\frac{2(p+1)}{(p-1)^2}w+|w|^{p-1}w
-\frac{(p+3)}{p-1}\partial_sw-2y\cdot \nabla\partial_s w,
\end{equation} 
\begin{equation}\label{defro0}
\mbox{where }\q L w = \frac 1\rho \dv(\rho \nabla w -\rho(y\cdot \nabla w)y)\mbox{ and }\rho(y) = (1-|y|^2)^\alpha\mbox{ with }\alpha = \frac 2{p-1}-\frac{N-1}2.
\end{equation}
By Antonini and Merle \cite{AMimrn01} in the subconformal range
\eqref{condpwave} and Merle and Zaag \cite{MZma05} in the conformal case
\eqref{condpwave'}, we have the following Lyapunov functional for
equation \eqref{eqw-ondes0} (obtained by multiplying the equation
by $\partial_s w \rho$ and then integrating in space):
\begin{equation}\label{defenergy0}
E_w(w)= \int_{|y|<1} \left(\frac 12 (\ps w)^2 + \frac 12  \left|\nabla w\right|^2 -(y\cdot \nabla w)^2+\frac{(p+1)}{(p-1)^2}w^2 - \frac 1{p+1} |w|^{p+1}\right)\rho dy.
\end{equation}
Note that this functional is defined in the unit ball $\{|y|<1\}$,
which corresponds to the backward light cone $\q C_{x_0, T(x_0),1}$ \eqref{defcone} in
the original variables $(x,t)$, thanks to the similarity variables transformation
\eqref{defwondes}. The functional $E_w(w)$ is decreasing
in the sense that for any $s\ge - \log T(x_0)$,
it holds that
\begin{align}
  \frac d{ds}E_w(w(s)) =&
- 2\alpha \int_{|y|<1}(\partial_s w(y,s))^2\frac{\rho(y)}{1-y^2} dy
  & \mbox{ for }p<p_c, \label{dissenergy}\\
  \frac d{ds}E_w(w(s)) = &
- \int_{|y|=1}(\partial_s w(\sigma,s))^2 d\sigma
  & \mbox{ for }p=p_c. \label{dissenergy'}                       
  \end{align}
  Note that $\alpha>0$ when $p<p_c$,
  and that the second identity is the limiting case of the first, when
  $p\to p_c$, in the distribution sense. Note that the dissipation
  degenerates to the boundary of the unit ball in the conformal case
  \eqref{dissenergy'}, which makes the analysis more difficult, in
  comparison with the subconformal case \eqref{dissenergy}. 
  Note also that our calculation in
  \eqref{dissenergy} breaks down in the superconformal range:
  \begin{equation}\label{supconf}
N\ge 2 \mbox{ and }p>p_c.
  \end{equation}
  In that case, assuming further that $p$ is Sobolev subcritical, Killip,
  Stovall and  Vi\c san \cite{KSVma14} derived some upper bounds on
  the blow-up rate. With Hamza, we improved that result in 
  \cite{HZdcds13} and \cite{HZjfa25}, with a different approach. In
  fact, we think all those bounds are not sharp.

  \medskip
  
Going back to the subconformal and conformal cases, note that the
existence of such a decreasing functional in the context of dispersive
equations like \eqref{eqnlw} is somehow surprising, given that
dispersive equations are time reversible. The solution of this
apparent paradox
resides in the fact that the time symmetry is lost after the
similarity variables transformation \eqref{defwondes}, which mixes up time and space,
leading to dissipative terms $\partial_s w$ and $\nabla \partial_s w$
in equation \eqref{eqw-ondes0}, unlike the original equation \eqref{eqnlw}.
The functional in \eqref{defenergy0} is in fact a crucial tool for the
blow-up analysis in \eqref{eqnlw}.
We would like to add that
dissipation in \eqref{dissenergy}-\eqref{dissenergy'} is the fundamental identity that
connects the blow-up study for this hyperbolic equation to the
framework we presented in Subsection \ref{secchaleur} in the parabolic
case (see \eqref{defE} and \eqref{dissipE} for the parabolic versions
of \eqref{defenergy0} and \eqref{dissenergy}-\eqref{dissenergy'}).

\medskip

Together with a blow-up criterion for equation \eqref{eqw-ondes0} from
\cite{AMimrn01}, the decay of the Lyapunov functional $E_w(w)$
introduced above allows to derive the following bound: 
\begin{equation}\label{boundE}
0 \le E_w(w(s)) \le C_0.
\end{equation}
Multiplying equation \eqref{eqw-ondes0} by $w \rho$ and then
integrating in space, we obtain a second energy identity. Combining
this identity with \eqref{boundE} and using various interpolation
inequalities (Sobolev, Hardy, Gagliardo-Nirenberg) allows us to
derive the desired bound on $w(y,s)$. 
Going back to the $u(x,t)$ setting, through the similarity variables
transformation \eqref{defwondes}, then using some covering technique
between light cones, we obtain the following statement in the
non-characteristic case, obtained in \cite{MZajm03}, \cite{MZma05} and
\cite{MZimrn05}: 
\begin{theorem}[Blow-up rate for $N\ge 1$ in the non-characteristic case] \label{theo1}
When $a\in \RR$ and $t\in [0, T(a))$, we have  
\[
 \dd\frac 1C
\le (T(a)-t)^{\frac 2{p-1}}\dd\frac{\|u(t)\|_{L^2(B(a, T(a)-t))}}{(T(a)-t)^{N/2}}
+(T(a)-t)^{\frac {p+1}{p-1}}
\left(
\dd\frac{\|{\partial_t u}(t)\|_{L^2(B(a, T(a)-t))}}{(T(a)-t)^{N/2}}+\dd\frac{\|\nabla u(t)\|_{L^2(B(a, T(a)-t))}}{(T(a)-t)^{N/2}}\right)\le C.
  \]
\end{theorem}
\noindent \textit{Remark}: The fractions appearing in this statement
are in fact $L^2$ averages of $u$, $\partial_t u$ and $\nabla u$ on
sections of the backward light cone $\q C_{a,T(a),1}$. The restriction
to that cone is compliant with the finite speed of propagation. The
$L^2$ norms are natural, given that the Cauchy problem for equation
\eqref{eqnlw} is solved in $H^1 \times L^2$.\\
\textit{Remark}: In the characteristic case, we proved in
\cite{MZimrn05} a similar though weaker result, where we have no lower
bound $\dd\frac 1C$, and where $L^2$ norms are taken on balls which are not
covering the whole section of the light cone, typically with radius
$\frac{T(a)-t}2$. In fact, as we will learn later from the asymptotic
behavior near characteristic points
when $N=1$ in item (ii) of Theorem \ref{threg}, and also for $N=2$ in
Section \ref{secpyr} for the case of the pyramid, the ``interesting''
objects for the asymptotic behavior in the characteristic case occur
near the boundary of the unit ball in the $w(y,s)$ setting, i.e. near
the boundary of the light cone in the $u(x,t)$ setting.\\
\textit{Remark}: We successfully generalized this result with Hamza in
\cite{HZnonl12} and \cite{HZjhde12} to perturbations of equation
\eqref{eqnlw} with lower order terms. Remarkably, we are able to prove a
similar statement for equations where the main term is not a pure
power nonlinearity (see \cite{HZna21} with Hamza and \cite{RZnodea25}
with Roy).

\subsection{Asymptotic behavior at blow-up.} \label{secasymp}

Take $N=1$ here. In higher dimensions, we were able to prove the same
results in the radial setting, assuming we are away from the origin
(see Merle and Zaag \cite{MZbsm11} and also Hamza and Zaag
\cite{HZnonl12} and \cite{HZjhde12} for some perturbative results).

\medskip

We first explain the key ideas behind our strategy. Once again, we
are inspired by the parabolic case. Indeed, since the Lyapunov
functional $E_w(w(s))$ \eqref{defenergy0}  is decreasing and bounded
from below (see \eqref{dissenergy} and \eqref{boundE}), it has a limit as $s\to \infty$.
In particular, from its dissipation given in \eqref{dissenergy}, it is
reasonable to think that $\partial_sw(y,s)\to 0$.
Thus, one expects $w(y,s)$ to converge to some stationary solution of equation
\eqref{eqw-ondes0}. It happens that when $N=1$, we can characterize
all stationary solutions in the energy space: $0$ and $\pm
\kappa(d,y)$, where 
\begin{equation}\label{defk}
  \forall |d|<1,\;\;\forall |y|<1,\;\;
\kappa(d,y)=\kappa_0 \frac{(1-|d|^2)^{\frac 1{p-1}}}{(1+d\cdot y)^{\frac 2{p-1}}},\mbox{ with }\kappa_0 = \left(\frac{2(p+1)}{(p-1)^2}\right)^{\frac 1{p-1}}
\end{equation}
(note that we give this definition for any $N\ge 1$ for the sake of
Subsection \ref{secpyr},
and that in the present case,
the inner product $d\cdot y$ simply becomes the real product $dy$). We
refer to such a solution as a ``soliton'', for two reasons which we
give in a remark after our next statement.  

\medskip

With this rough explanation, we are able to assert that our techniques
allow us to 
    completely prove the \textit{``Soliton
      Resolution Conjecture''} for equation \eqref{eqnlw} in one
    space dimension, in the sense that:\\
    - near a non-characteristic point, the solution approaches one
    soliton;\\
    - near a characteristic point, the solution decomposes into a sum
    of decoupled solitons.

\medskip

More precisely, this is our statement from  \cite{MZjfa07}, \cite{MZajm12} and
\cite{CZcpam13}:   
\begin{theorem}[Asymptotic behavior near the blow-up
  graph]\label{threg} 
  Take $N=1$ and consider
  $a\in \m R$.\\
  (i) {\bf Case where $\xx\in\RR$: Convergence to a soliton}.
  It holds that
  \begin{equation}\label{profile}
w_a(y,s) \to \pm \kappa(d(a),y)\mbox{ as } s\to \infty,
  \end{equation}
  for some $|d(a)|<1$. \\
  (ii) {\bf Case where $\xx\in\SS$: Decomposition as a multi-soliton}.
  It holds that
\begin{equation}\label{cprofile00}
  w_a(y,s) \sim \pm \dd\sum_{i=1}^{k(\xx)} (-1)^{i+1}\kappa(d_i(s),y)
  \mbox{ as } s\to \infty,
\end{equation}
for some integer $k(a) \ge 2$, where
$d_i(s)=-\tanh \zeta_i(s)\in (-1,1)$ and
\begin{equation}\label{equid}
\zeta_i(s)=\left(i-\frac{(k(\xx)+1)}2\right)\frac{(p-1)}2\log s + C_i,
\end{equation}
for some explicit constants $C_i$.
\end{theorem}
\textit{Remark}: The convergence in \eqref{profile} holds in
$H^1(-1,1)$, hence in $L^\infty(-1,1)$ from the Sobolev embedding. As
for \eqref{cprofile00}, it holds in the space $\q H_0$ such that
\begin{equation}\label{defnh}
\|q\|_{{\q H}_0}^2\equiv \int_{|y|<1} \left(q^2+|\nabla q|^2-(y\cdot
  \nabla q)^2)\right)\rho dy
\end{equation}
(again, we give this definition for all $N\ge 1$ for further use in
Subsection \ref{secpyr}).\\
\textit{Remark}: The constants $C_i$ in item (ii) are uniquely determined (up to some
additive constant) by the fact that the following system holds:
\begin{equation}\label{toda1}
  \forall i=1,\dots,k,\;\;
  \frac 1{A} \zeta_i'(s) = e^{-\frac
    2{p-1}(\zeta_i(s)-\zeta_{i-1}(s))}-e^{-\frac
    2{p-1}(\zeta_{i+1}(s)-\zeta_i(s))},
\end{equation}
where $A>0$, $\zeta_0(s) \equiv - \infty$ and $\zeta_{k+1}(s) \equiv + \infty$.  
This system was first exhibited with Merle in \cite{MZajm12}. For the
explicit value of $A$, see Azaiez, Jendrej and Zaag \cite{AJZ25}.
We refer to this system as
a  ``first order Toda system'', given the similarity with the Toda lattice system,
where the first time derivative is replaced by a second derivative.\\
\textit{Remark}: This statement remains valid for the complex-valued
case of equation \eqref{eqnlw}, where $u\in\m C$, up to replacing
the $\pm$ sign in items (i) and (ii) by $e^{i \theta_0}$ for some
$\theta_0 \in \m R$ (see Azaiez \cite{Atams15},
Azaiez, Jendrej and Zaag \cite{AJZ25}). In other words, the blow-up
behavior is asymptotically scalar. In addition, statement (i) is valid
also in the vector-valued case $u\in \m R^m$ with $m\ge 3$ (see Azaiez and Zaag \cite{AZbsm17}).\\ 
\textit{Remark}: We adopt the terminology ``soliton'' for
$\kappa(d,y)$ \eqref{defk} for two reasons linked to the following
change of variables:
\[
  \bar w_a(\xi,s)= (1-y^2)^{\frac 1{p-1}}w_a(y,s)
  \mbox{  with }y=\tanh \xi. 
\]
Indeed, in this setting, estimate \eqref{cprofile00} above translates
as 
\begin{equation}\label{resultinxi}
\|\bar w_a(\xi,s)-
\pm\kappa_0\sum_{i=1}^{k(s)}(-1)^{i+1} \cosh^{-\frac 2{p-1}}(\xi-\zeta_i(s))\|_{H^1\cap L^\infty(\m R)} 
\to 0 \mbox{ as }s\to \infty,
\end{equation}
where we see that the $\kappa(d_i(s),y)$ become space shifts of the
bump $\kappa_0 \cosh^{-\frac 2{p-1}}\xi$, centered at $\zeta_i(s)$ given above
in \eqref{equid}. Since the distance between the bumps is given
by $\frac{(p-1)}2 \log s\to +\infty$, we see that the bumps are
``decoupled''. In particular, they are evolving ``alone'' or
``isolated'', which justifies to call them solitons. In addition, the
bumps in the $\xi$ setting satisfy the equation of the solitons of the
Korteweg-de Vries (KdV) equation, which gives a second reason to call
them ``solitons'' (see \cite{KdVpm95}).\\
\textit{Remark}: In the characteristic case, we see that $\zeta_i(s)
\to \pm \infty$, hence $d_i(s) \to
\pm 1$, except maybe when $i = \frac{k+1}2$ and $k$ is odd. In view of
the decomposition \eqref{resultinxi} in the $\bar w(\xi,s)$ setting,
using the similarity variables transformation \eqref{defwondes}, this
means that the solitons are escaping towards the boundary of the light
cone $\q C_{a, T(a),1}$.

\medskip

When $N=1$, as a consequence of the blow-up behavior in similarity
variables (see Theorem \ref{threg} below),
we have the following result on the blow-up rate (in $L^\infty$) in
sections of the backward light cone with vertex $(a, T(a))$:
\begin{corollary}[Blow-up rate for $N=1$]\label{corspeed}\\
(i) For all $a\in\RR$ and $t\in[0, T(a))$, we have
\[
\frac 1{C(T(a)-t)^{\frac 2{p-1}}} \le \sup_{|x-a|<T(a)-t}|u(x,t)|\le \frac C{(T(a)-t)^{\frac 2{p-1}}}.
\]
(ii) For all $a\in\SS$ and
$t\in[0, T(a))$, we have:
\[
\frac{|\log(T(a)-t)|^{\frac{k(a)-1}2}}{C(T(a)-t)^{\frac 2{p-1}}}\le \sup_{|x-a|<T(a)-t}|u(x,t)|\le \frac{C |\log(T(a)-t)|^{\frac{k(a)-1}2}}{(T(a)-t)^{\frac 2{p-1}}},
\]
where the integer $k(a)\ge 2$ corresponds to the number of solitons
appearing in the decomposition of $w_a(y,s)$ given above in
\eqref{cprofile00}. 
\end{corollary}
\noindent {\it Remark}: Note that the blow-up rate in the
non-characteristic case is given by the solution of the associated ODE
$u''=u^p$, whereas it is higher in the characteristic case. 

\subsection{Regularity of the blow-up graph.}\label{secreg}

We address this issue when $N=1$.
As we will see in a moment, the regularity of the blow-up graph is a
consequence of the asymptotic behavior given in Subsection
\ref{secasymp}
above. For that reason, we will handle the case of non-characteristic
and characteristic points separately.

\medskip

If $a\in \RR$, using Theorem \ref{threg} and going back to the
$u(x,t)$ setting, we write from
\eqref{defwondes}
\begin{equation}\label{uprofile}
  u(x,t) \sim \pm \kappa_0 \frac{(T(a)-t)^{-\frac 1{p-1}}}
  {(T(a) - t+d(a)(x-a))^{\frac 2{p-1}}} \mbox{ as } t\to T(a),
  \end{equation}
  uniformly in sections of the backward light cone $\q C_{a,T(a),1}$
  \eqref{defcone}. It happens that the function on the right-hand side
  of \eqref{uprofile} is a particular solution of equation
  \eqref{eqnlw} whose blow-up graph is the line of equation
  \[
t= T(a) +d(a)(x-a),
  \]
  which has $d(a)$ as a slope. Because of \eqref{uprofile}, the
  blow-up graphs of both solutions have to be tangent at $a$, hence
  $T'(a) = d(a)$. The stability with respect to $a$ of $\kappa(d,y)$ as a
  blow-up profile for $w_a$ is a crucial ingredient. It is proved in
  \cite{MZcmp08} (see also \cite{MZcmp15} for $N\ge 2$). This
  stability implies that $\RR$ is open.

  \medskip

  When $a\in \SS$, consider $b\in \m R$ close to $a$. Our idea is to show
  that unlike $w_a$, $w_b$ will not decompose as a multi-soliton for
  $s\to \infty$. It will instead converge to a single soliton, meaning
  that we fall in item (i), which implies that $b$ is
  non-characteristic, hence, $a$ is an \textit{isolated}
  characteristic point. 

  \medskip

The starting point is the following direct connection between $w_b$
and $w_a$, which follows from a double application of the similarity
variables transformation \eqref{defwondes}:
\begin{equation}\label{wawb}
    w_b(z,\sigma_0) =e^{\frac 2{p-1}(s_0-\sigma_0)} w_a(y,s_0),
    \mbox{ with } b+z e^{-\sigma_0} = a +y e^{-s_0} 
\mbox{ and } T(b) - e^{-\sigma_0} = T(a) - e^{-s_0}.
  \end{equation}
  Note that $(z,\sigma_0) \to (y,s_0)$ as $b\to a$.

  \medskip

  Since we know from item (ii) in Theorem \ref{threg} that
  $w_a(y,s_0)$ is close to a sum of decoupled solitons
  $\pm\kappa(d_i(s_0),y)$,
  identity \eqref{wawb}
  also shows that $w_b(z,\sigma_0)$ is a sum of transformations of
  those solitons by \eqref{wawb}, which happen to be members of the
  following family of ``generalized'' solitons:
  \[
\kappa^*(d,\nu,y) = 
\kappa_0 \frac{(1-d^2)^{\frac 1{p-1}}}{(1+\nu+dy)^{\frac 2{p-1}}},
  \]
  When $\nu = \mu e^s$, $\kappa^*(d,\mu e^s, y)$ is in fact a
  particular solution of equation \eqref{eqw-ondes0}, with the
  following property:\\
  - if $\mu>0$, then $\kappa(d,y)>\kappa^*(d,\mu e^s,y) \to 0$ as
  $s\to \infty$ ;\\
  - if $\mu<0$, then $\kappa(d,y) <\kappa^*(d,\mu e^s,y)$ which grows
  and blows up in some finite time $S^*=- \log\left(\frac{|d|-1}\mu\right)$.

  \medskip

  From symmetry, we may take $b<a$.
  It happens that $w_b(z,\sigma_0)$ is a sum of
  $\pm \kappa^*(\bar d_i(s_0), \nu_{0,i}, z)$, with $\nu_{0,i} \ge 0$, except maybe for
  $i=1$. In other words, all the solitons, except number $i=1$
  are smaller than the regular soliton
  $\pm \kappa(\bar d_i(s_0), z)$, which indicates they will decay then
  vanish for
  $s\ge \sigma_0$.
  At the end,
only the soliton number $i=1$ remains, 
which makes us fall in iten (i) of Theorem
  \ref{threg}, namely in the non-characteristic case. This
  ``losing-solitons'' strategy is performed in \cite{MZdmj12}. It
  relies on a good understanding of the dynamics of equation
  \eqref{eqw-ondes0} near the multi-soliton, in particular its instability.

  \medskip

\begin{figure}[htbp]\label{corner-shaped}
  \centering
  \includegraphics[width=0.5\columnwidth]{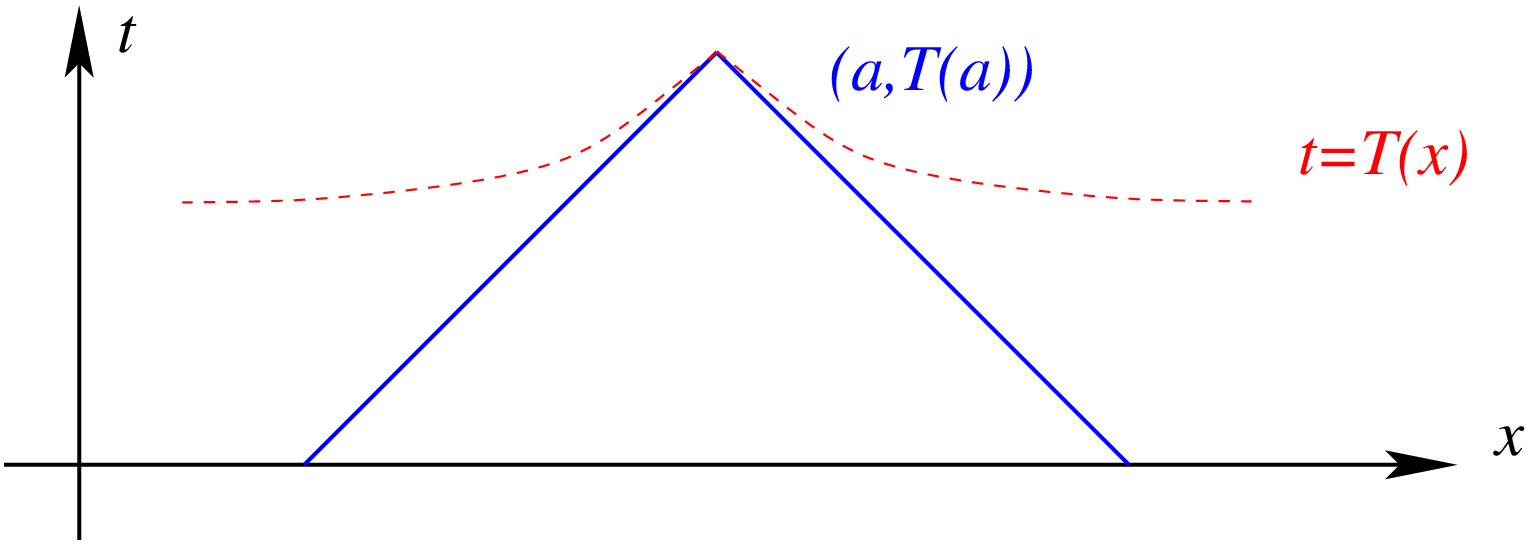}
  \caption{The blow-up set is corner-shaped near characteristic points.}
\end{figure}

  In the following statement, we summarize our results from \cite{MZcmp08}, \cite{MZajm12} and \cite{MZdmj12}:
\begin{theorem}[Geometry of the blow-up graph]\label{old} $ $\\
(i) The set $\RR$ is non-empty and open, and $x\mapsto T(x)$ is of
class $C^1$ on $\RR$. Moreover, $\forall a\in \RR$, $T'(a) = d(a)$
shown in the asymptotic behavior of $w_a$ given in \eqref{profile}.\\
(ii) The set $\SS$ is made of isolated points, and given $a\in\SS$, if $0<|x-a|\le \delta_0$, then
\begin{equation}\label{chapeau0}
\dd\frac{|x-a|}{C_0|\log(x-a)|^{\frac{(k(\xx)-1)(p-1)}2}}\le T(x)- T(a)+|x-a| \le \frac{C_0|x-a|}{|\log(x-a)|^{\frac{(k(a)-1)(p-1)}2}}
\end{equation}
for some $\delta_0>0$ and $C_0>0$, where $k(a)\ge 2$ is the number of
solitons appearing in the asymptotic decomposition of the solution in
\eqref{cprofile00}. In particular, $T(x)$ is right and left
differentiable at $a$, with $T'_l(a)=1$ and $T'_r(a)=-1$ (see Figure \ref{corner-shaped}).
\end{theorem}
\noindent{\it Remark}: In \cite{Ncpde08}, Nouaili improves the regularity of the restriction of $x\mapsto T(x)$ to $\RR$ to $C^{1,\alpha}$ for some $\alpha>0$.\\
{\it Remark}: The fact that $\SS$ is made of isolated points certainly
does not hold in general for nonlinear wave equations. A
counter-example is given by Alinhac in  \cite{Anum06} for
this equation
$\partial^2_{t} u =\partial^2_{x} u+\partial_x u\partial_t u$. 

\subsection{Construction of examples and existence results related to
  characteristic points.} \label{seccons}

Again, we consider $N=1$ here. 
Note first that it was commonly conjectured before our contributions that $\SS$ was empty. In particular, the first result in this direction is due to Caffarelli and Friedman in \cite{CFtams86} and \cite{CFarma85}. Using the maximum principle, they show that $\SS = \emptyset$, under some restrictive conditions on initial data that ensure the fact that
\[
u\ge -C\mbox{ and }
\partial_t u \ge (1+\delta_0)|\partial_x u|\mbox{ for some }C>0\mbox{ and }\delta_0>0.
\]
In \cite{MZajm12}, we derived the first example of a solution with
$\SS\neq \emptyset$:
\begin{proposition}[Existence of initial data with $\SS\neq \emptyset$]
   \label{propelem}
When $N=1$, if the initial data $(u_0,u_1)$ is odd with larger
plateaus where it is equal to some blow-up solution of the associated ODE
$u''=u^p$, then, $u(x,t)$ blows up in finite time, $u(0,t)=0$, and the origin is a
characteristic point.  
\end{proposition}

Following this and the classification result given in Theorem
\ref{threg}, a natural question arises, in accordance with the duality
between the classification and the  construction questions we
mentioned in Subsection \ref{secquest}: 
\medskip

\textit{Do we have examples of initial data for every blow-up modality given in Theorem \ref{threg}?} 

\medskip

The answer is yes, as we proved in C\^ote and Zaag \cite{CZcpam13}: 
\begin{theorem}[Construction of examples for all the blow-up
  modalities in Theorem \ref{threg}]$ $\\
  (i) For any $a\in \m R$ and $d\in(-1,1)$, there exists a solution $u(x,t)$ to
equation \eqref{eqnlw} which blows up in finite time such that $a\in
\RR$ and the convergence in \eqref{profile} holds for $w_a$ with
$d(a) = d$.\\
(ii) Given an integer $k\ge 2$ and $a\in \m R$, there exists a solution $u(x,t)$ to
equation \eqref{eqnlw} which blows up in finite time such that $a\in
\SS$ and the multi-soliton decomposition in \eqref{cprofile00} holds for $w_a$ with
$k(a) =k$.
\end{theorem}
\textit{Remark}:
The proof of item (i) is trivial, thanks to the finite speed of
propagation and to the fact that $\kappa(d,y)$ is a particular solution for equation
\eqref{eqw-ondes0}.
As
for item (ii), the proof relies on our good understanding in
\cite{MZdmj12} of the dynamics of equation \eqref{eqw-ondes0} near
multi-solitons.\\
{\it Remark}: Some refinements for item (ii) are given in \cite{HZjde19} in
collaboration with Hamza.\\
\textit{Remark}: Many papers have addressed the question of
multi-solitons for various PDEs. See for example
Jendrej and Martel \cite{JMjmpa20},
Martel and Merle \cite{MMarma16}.

\medskip

\subsection{Nearly Pyramidal blow-up graph in higher dimensions.} \label{secpyr}
We take $N\ge 2$ in this subsection.
Blow-up results in
this case
are not that advanced as for
$N=1$. Indeed, the blow-up rate near non-characteristic points is the
only available result for general initial data (see Theorem \ref{theo1}).
For the asymptotic behavior and the regularity of the blow-up graph, 
the only known result is at non-characteristic points, where we show
in \cite{MZtams16} and \cite{MZcmp15} that
the blow-up graph
is $C^1$,
 under a reasonable assumption on the profile, namely that
the local blow-up profile is given by $\kappa(d,y)$ \eqref{defk}.
The radial case outside the origin is also completely understood in \cite{MZbsm11}, since it reduces to a perturbation of the one-dimensional case.\\
Concerning the behavior of radial solutions at the origin, Donninger and Sch{\"o}rkhuber were able to prove the stability of the space-independent solution 
with respect to perturbations in initial data, in the Sobolev subcritical range \cite{DSdpde12} and also in the supercritical range in \cite{DStams14}. Some numerical results are available in a series of papers by Bizo\'n and co-authors (see \cite{Bjnmp01}, \cite{BCTnonl04}, \cite{BZnonl09}).  
See also Killip and Vi\c san \cite{KV11}.

\medskip

Our main focus in this subsection concerns 
the existence of blow-up
solutions to equation \eqref{eqnlw} with $\SS\neq \emptyset$.
As asserted above in Subsection \ref{seccons}, we have examples in one
space dimension with $\SS\neq \emptyset$. 
Those examples naturally extend to the radial case outside the origin,
as we showed in \cite{MZbsm11}. Thanks to the 1-d and radial cases, we
naturally derive trivial $2$-dimensional solutions, where $\SS$ is either a line, or a circle.
From the finite speed of propagation, we may have parallel lines or concentric circles, and the local blow-up behavior is always rigorously one dimensional.
In particular, no example is known in higher dimensions, with $\SS$
locally reduced to an isolated point. In \cite{MZ16}, we give such an example.

\subsubsection{Statement of the results for $N=2$.}

Before stating our result, let us introduce
\begin{equation*}
\bar d(s) = - \tanh \bar \zeta(s)\mbox{ where }\bar \zeta(s) = \left(\frac{p-1}4\right) \log s -\frac{(p-1)}4 \log \left(\frac{p-1}{4\barc}\right)
\end{equation*}
which is an explicit solution to the ODE
\[
\frac 1{\barc} \frac{d\bar\zeta}{ds} = e^{-\frac 4{p-1}\bar \zeta}
\]
where $A=A(p)$ is the same constant as in system \eqref{toda1}
relevant for multi-solitons in 1-d.
Note that we have for some $C_0(p)>0$,
\begin{equation}\label{devdbar}
  1+\bar d(s)
  \sim C_0s^{-\frac{p-1}2}\mbox{ as }s\to \infty.
\end{equation}
 
\medskip
Let
$(e_1, e_2)$ 
be
the canonical basis of $\m R^2$.
This is the statement of our result:
\begin{theorem}[Existence of a blow-up solution with 
an isolated characteristic blow-up point and a blow-up surface which
is nearly pyramidal]
\label{mainth} There exists $u(x,t)$ a solution to equation \eqref{eqnlw},
which is symmetric with respect to the axes and anti-symmetric with
respect to the lines $\{y = \pm  x\}$,
with the following properties:\\ 
{\bf  (A) (Blow-up with an isolated characteristic point)}. The
solution $u(x,t)$ blows up on some blow-up graph $\Gamma =
\{(x,T(x)\}$, and for some $\delta>0$, we have $\SS\cap B(0, \delta) =
\{0\}$.\\ 
 {\bf (B) (The blow-up graph is nearly a pyramid)}. The function $T$
 is symmetric with respect to the axes and the
 lines $\{y = \pm  x\}$,
 and $C^1$ outside the
 lines $\{y = \pm  x\}$,
Moreover, when $0\le x_2<x_1\le\delta$, we have for some $C_0=C_0(p)>0$:
\[
  T(x)=T(0)-x_1(1-C_0|\log x_1|^{-\frac{p-1}2})+o(x_1|\log x_1|^{-\frac{p-1}2})+o(x_2|\log x_1|^{-\frac{p-1}4}).
\]
{\bf (C) (Blow-up behavior of the solution)}.
We have the following behavior for $w_x$ for $0\le x_2\le x_1\le \delta$
as $s\to \infty$:\\
(i) if $x=0$, then 
\begin{equation}\label{cprofile0}
\left\|w_{0}(y,s)-\left(\kappa(\bar d(s)e_1,y)+\kappa(-\bar d(s)e_1,y) - \kappa(\bar d(s)e_2,y) - \kappa(-\bar d(s)e_2,y)\right)\right\|_{\q H}\to 0
\end{equation}
where $1+\bar d(s)\sim C_0s^{-\frac{p-1}2}$ as $s\to \infty$;\\
(ii) if 
$x_2<x_1$, 
then $w_x(s)$ converges as $s\to +\infty$ to 
$\kappa(d(x)e_1)$, 
with
\begin{equation*}
d(x)+1\sim C_0|\log x_1|^{-\frac{p-1}2}\mbox{ as }x\to 0. 
\end{equation*}
(iii)  if $x\neq 0$ with $x_1= x_2$, then $w_x(s_n)$ converges to some
stationary solution $w^*_x$ for some sequence $s_n \to \infty$, where
$w^*_x$ is a genuinely two-dimensional stationary solution of equation
\eqref{eqw-ondes0}. 
\end{theorem}
\textit{Remark}: 
The existence of the new stationary solution of equation
\eqref{eqw-ondes0} just mentioned at the end of this theorem follows
from an indirect argument we use when $x$ is on
the lines $\{ y = \pm x\}$.\\
\textit{Remark}: 
Note 
from the symmetries of the solution that we
have $u(x,t)=0$ on
the lines $\{ y = \pm x\}$.
In one space dimension, such a property implies that $x$ is a
characteristic point. Surprisingly, in our two-dimensional setting,
only the origin is a characteristic point, and the other points on
the lines $\{ y = \pm x\}$ are non-characteristic points, showing a genuinely two-dimensional behavior.

\subsubsection{Generalization and extensions of the result.}
 Our result can be generalized to other pyramids, with any regular polygon as a basis. 
In higher space dimensions $N\ge 3$, we naturally generalize our results to a pyramid with a hypercube as a basis. Moreover, 
using a Lorentz transform near  some point of
the lines $\{ y = \pm x\}$
different from the origin, we can tilt the blow-up surface and obtain the existence of a blow-up solution of equation \eqref{eqnlw}, with a tent-shaped (at the first order) blow-up surface, no characteristic point in some neighborhood, a slope approaching $\frac{\sqrt 2}2$ and an upper edge depending on $x_2$. This tent is in fact new and different from the one obtained by considering a solution depending only on $x_1$ with a characteristic point at the origin. Indeed, in two space dimensions, this ``naive'' tent has a line of characteristic points on its upper edge, a slope approaching $1$, and an upper edge that does not depend on $x_2$.
\subsubsection{The strategy of the proof.}

Our proof relies on 3 main steps, which we present in the following
paragraphs (for details and proofs, see \cite{MZ16}).

\medskip

\paragraph{Step 1: Construction of a solution for equation \eqref{eqnlw} showing 4 solitons in the backward light cone.}

In this step, 
we construct a blow-up solution to equation \eqref{eqnlw} defined only in the backward light cone with vertex $(0,T(0))$ and showing 4 solitons for $w_0$ at the origin as in \eqref{cprofile0}. 
Then, using the finite speed of propagation, we derive from the latter
a blow-up solution to the Cauchy problem of equation
\eqref{eqnlw}. Note that this construction step follows the
classical scheme of a {\it ``construction with a prescribed
  behavior''}, which proved to be efficient for various PDEs (see the
introduction of \cite{MZcpam18}), in particular in our construction of
a multi-soliton for the same equation in 1-d with C\^ote in \cite{CZcpam13}.

\medskip

\paragraph{Step 2: Instability of the 4-soliton solution when the vertex of the backward light cone leaves the origin.}
This step 
is the very heart of our argument. Here, we aim at understanding the
instability of the 4-soliton solution of equation \eqref{eqw-ondes0}
we have for $w_0$ \eqref{cprofile0}, when we move outside the origin
to consider the behavior of $w_{x_0}$ where $x_0\neq 0$. 
We did that already in one space dimension in \cite{MZdmj12}, as we
explained above in Subsection \ref{secreg}. In 2-d, the strategy is
basically the same. However, the situation is much more delicate than
in 1-d, mainly because of the dynamics at the lines $\{x_{0,1} = \pm x_{0,2}\}$.

\medskip

\paragraph{Step 3: The local behavior of $T(x_0)$ in connection with the dynamics of $w_{x_0}$.}

As usual with blow-up problems (heat, wave), the asymptotic behavior of the solution at
blow-up and the regularity of the blow-up set are linked and advanced side by side in the
proof (see \cite{Zihp02}, \cite{Zcmp02}, \cite{Zdmj06} and \cite{Zbeit00} for the semilinear heat equation; see \cite{MZjfa07}, \cite{MZcmp08}, \cite{MZajm12}, \cite{MZdmj12},  \cite{MZxedp10} and \cite{MZsnp09} for the semilinear wave equation).

The present situation is no exception. As a matter of fact, in this step, 
we make the link between the dynamics of $w_{x_0}$ from the previous
section and the local behavior of $T(x_0)$. This is in fact the new
feature of our paper \cite{MZcpam18}, which makes it original. In
particular, we use families of moving non-characteristic cones
together with subtle elementary geometric methods to derive the nearly
pyramidal shape of the blow-up surface. The delicate case is the case
where $x_0$ is on
the lines $\{ y= \pm x\}$ 
since this is a new situation, not encountered in dimension 1. It is worth noticing that the moving cone technique simplifies the ``moving plane'' technique we use earlier in one space dimension in \cite{MZcmp08}.

\section{Perspectives and open problems.}

As we have already mentioned in the introduction, our focus on
singularity formation for (NLW) stemmed from a broader interest in
more physical problems, with less ``nice'' properties, or less
``structure''.

\medskip

Our aim was to develop new tools for (NLW), with the hope to apply
them to other equations showing blow-up. Our analysis was largely
successful in one space dimension. However, many open problems still
remain for (NLW). 

\medskip

Accordingly, two possible directions remain to be explored, following
this presentation:

- Advancing knowledge about blow-up for (NLW);

- Extending the techniques to other models, with less structure
and more physical interest.

\medskip

We will further comment on these directions in the following.

\subsection{Future directions for (NLW).}

The situation is very much dependent on the exponent $p$ in
\eqref{eqnlw}.

\medskip

In the subconformal range \eqref{condpwave}, the situation remains
largely unclear for $N\ge 2$. Indeed, the main obstruction lays in the
absence of a classification of \textit{all} finite energy solutions to the stationary
problem of equation \eqref{eqw-ondes0}, which happens to be a
degenerate elliptic problem posed in the unit ball. We already know some
solutions, namely the solitons $\kappa(d,y)$ given in \eqref{defk},
and also the countable family of radial solutions given by Biz\'on,
Breitenlohner, Maison and Wasserman in \cite{BBMWnonl10} and
\cite{BMWnonl07}, when $N=3$, $p=3$ or $p\ge 7$ and odd. Classifying
all the solutions should open the road to a complete classification of
the blow-up behavior when $p<p_c$ and $N\ge 2$.

\medskip

The conformal case $p=p_c$ should follow the same pattern, though we
expect some complications coming from the degeneracy of this case, as
one can already see from the dissipation of the Lyapunov functional in
\eqref{dissenergy'}, in comparison with the subconformal case in
\eqref{dissenergy}. In our previous papers, we encountered some
serious technical problems, preventing us from extending the results
of the subconformal case to the conformal case, in particular in
\cite{MZajm03} and \cite{MZtams16}. Some difficulties were overcome in
\cite{MZma05}, and many remain.

\medskip

In the superconformal and Sobolev critical range $p_c<p<p_S$, the
blow-up rate is already an open issue (see our comments right after
\eqref{supconf}.

\medskip

The Sobolev critical and supercritical case $p\ge p_S$ go
beyond the scope of this presentation, and many contributions already
exist; some of them were already cited in this presentation. Nevertheless, I would like to
cite Klainerman \cite{Kgafa00}, who mentioned the determination
of \textit{generic} behavior as an important project. I do think this
is still largely open, especially for $p>p_S$.

\subsection{Extension to more physically relevant models.}
As already mentioned, (NLW) shares important features with more
physically relevant models, in particular with General Relativity (GR)
(see Donninger, Schlag and Soffer \cite{DSScmp12}). On top of them, we
have the finite speed of propagation, which justified for (NLW) the
restriction of the analysis to backward light cones with vertex on the
blow-up graph. Such ideas should be important (though not enough) in GR. 

\medskip

A different challenge lays in the following version of the nonlinear
wave equation
\begin{equation}\label{eqeft}
  \partial_t^2 u = \partial_x^2 u +(\partial_x u)^2,
\end{equation}
which stems from effective field theories (EFT) originated from cosmological studies (see our paper
\cite{EHZnonl23} with Eckmann and Hassani). The main difference with (NLW) lays in the absence
of a variational structure for this equation, preventing the existence
of a conserved energy or Lyapunov functional. As the existence of the
Lyapunov functional \eqref{defenergy0} for the similarity variables
version \eqref{eqw-ondes0} was crucial in our analysis of blow-up for
(NLW), it clearly appears that our strategy breaks down for
\eqref{eqeft}. Nevertheless, we were able to give some blow-up results
in \cite{EHZnonl23}. We would like also to mention the paper by Ghoul,
Liu and Masmoudi in \cite{GLM25}.

In spite of those important results, we think that addressing blow-up
for equation \eqref{eqeft} along with other non variational hyperbolic
PDEs (such as the hyperbolic MEMS model \cite{DZm3as19}) should be an important research direction.

\section*{Acknowledgments.}
 The author would like to thank Frank Merle for his guidance throughout
my career, and for precious advice in writing this paper. Many thanks
also to Thomas Duyckaerts, Jean-Pierre Eckmann, Yvan Martel and Nejla Nouaili for
their valuable comments on the manuscript, and also to Irfan Glogi\'c
for pointing out some important references. The author also wishes to
thank Pierre Rapha\"el for his valuable support through his ERC
project SWAT. 

\bibliographystyle{siamplain}
\bibliography{../ref}
\end{document}